\documentclass[reqno,b5paper]{amsart}
\usepackage{amsmath,amssymb,amsbsy,amsfonts,amsthm,latexsym,amsopn,amstext,amsxtra,euscript,amscd,amsthm,graphicx}

\newtheorem{theorem}{Theorem}
\newtheorem{corollary}{Corollary}

\newcommand {\N} {\mathcal{N}}
\newcommand {\CC} {{\mathbb C}}

\newcommand{\C}{\mathbb C}

\newcommand{\acr}{\newline\indent}

\begin{document}

\title[30 years of collaboration]{30 years of collaboration\\ \footnotesize{\dots between the Austrian Diophantine Number Theory research group and the Number Theory and Cryptography School of Debrecen}}

\author[C. Fuchs \and L. Hajdu]{Clemens Fuchs* \and Lajos Hajdu**}

\address{\llap{*\,}University of Salzburg \acr Department of Mathematics\acr Hellbrunner Str. 34/I \acr 5020 Salzburg \acr AUSTRIA}
\email{clemens.fuchs@sbg.ac.at}
\address{\llap{**\,}University of Debrecen\acr Institute of Mathematics\acr P.O. Box 12\acr 4010 Debrecen, HUNGARY}
\email{hajdul@science.unideb.hu}

\thanks{C.F. was supported by FWF (Austrian Science Fund) grant No. P24574. L.H. was supported by the OTKA grants K100339 and K115479.}

\begin{abstract}
We highlight some of the most important cornerstones of the long standing and very fruitful collaboration of the Austrian Diophantine Number Theory research group and the Number Theory and Cryptography School of Debrecen. However, we do not plan to be complete in any sense but give some interesting data and selected results that we find particularly nice. At the end we focus on two topics in more details, namely a problem that origins from a conjecture of R\'enyi and Erd\H{o}s (on the number of terms of the square of a polynomial) and another one that origins from a question of Zelinsky (on the unit sum number problem), which will be presented in turn. This paper evolved from a plenary invited talk that the authors gave at the Joint Austrian-Hungarian Mathematical Conference 2015, August 25-27, 2015 in Gy\H{o}r (Hungary).
\end{abstract}

\maketitle

\section{Introduction}

There are a lot of ongoing scientific collaborations between Austrian and Hungarian mathematicians in various fields of mathematics. In this paper we shall give some account on this activity. However, since we can give only a very restricted view into the fortunately extremely rich and vivid research activity among Austrians and Hungarians, we shall restrict ourselves to research done in Number Theory. Further, we do not even claim at all that we are able to give a comprehensive account of the research done and/or in progress inside this traditional research field. Therefore, we shall concentrate on Diophantine problems inside Number Theory, and even more restrictive, we shall focus on the research done between the members of the research groups (in a broader sense) to which the authors belong to. What we hope to be able to do, is to try to exhibit some of the main lines describing the very fruitful and long standing cooperation of these research groups. Certainly, this is an enormous task in itself, surely we shall miss beautiful, important and interesting points, which of course is completely due to our own ignorance and for which we apologize to the colleagues not properly mentioned. Let us see what we can do about it!\\

The paper is organized as follows: In the next section we shall explain who ``we'' are - that is, we introduce our research groups. In Section 3 we shall give some interesting data on our collaboration. A closer look at some of the most important joint research topics and results within the collaboration will be given in four subsections in Section 3. In Sections 4 and 5, respectively, we shall make a closer look at two specific problems in some more details. The first of these topics comes from the more recent research of the first author, the second topic was chosen by the second author and it contains contributions of him with Austrian colleagues. In the last section, we will summarize the main points and give some outlook on current and future research. The references are split up in two parts: \cite{MR2075644} up to \cite{MR3327847} lists the references which belong to our collaboration, from \cite{MR1871104} onwards we list other cited references (in both blocks the references are listed in the usual alphabetical ordering).

\section{Who we are - the research groups we belong to}

The collaboration we are referring to is certainly a collaboration between Austrian and Hungarian mathematicians. As already mentioned in the introduction, we are more specific: We will restrict to joint work inside number theory. Even more specifically we focus on Diophantine number theory. Even more restrictively we focus on the collaboration between the two research groups (in a broader sense) to which we belong to. First we will introduce the Austrian group involved in the collaboration.\\

For the purpose of this paper we shall focus on the research groups in {\bf Graz} and {\bf Leoben} and their ``spin-offs''\footnote{This is meant in the sense of the people from our collaboration that we have in mind. C. Heuberger from the University of Klagenfurt and the first author from the University of Salzburg both did their PhD with R. Tichy at TU Graz.} in {\bf Klagenfurt} and {\bf Salzburg}. These are the groups in which the first author has his scientific origins (he was PhD student, assistant and postdoc in Graz as well as teaching assistant in Leoben, now he has his own research group in Salzburg). The members include (given in order receiving their PhD degree)\footnote{In this list we are not mentioning our PhD students.}:
\begin{itemize}
\item {\bf TU Graz}: R. Tichy, P. Grabner, S. Frisch, C. Elsholtz, C. Aistleitner\footnote{C. Aistleitner did his PhD with R. Tichy at TU Graz. He is currently in Linz, but was very recently awarded a START grant of the FWF to take place at TU Graz.}, C. Frei, D. Kreso, \dots, G. W\"ustholz\footnote{G. W\"ustholz is Emeritus Professor at ETH Z\"urich and is currently involved in the research project P26114 at TU Graz granted by the FWF.}, \dots
\item {\bf KFU Graz}: F. Halter-Koch, G. Lettl, A. Geroldinger, \dots
\item {\bf Leoben}: P. Kirschenhofer, J. Thuswaldner, C. van de Woestijne, \dots
\item {\bf Klagenfurt}: C. Heuberger\footnote{C. Heuberger has shifted his research interests slightly away from Diophantine number theory, but since he had an important role in the collaboration, we want to mention him nevertheless.}, \dots
\item {\bf Salzburg}: C. Fuchs, V. Ziegler, I. Pink, R. Paulin\footnote{R. Paulin and I. Pink are Hungarians who currently have postdoc positions (FWF grant P24574 resp. P24801) at the University of Salzburg.}, \dots, P. Hellekalek, W. Schmid, \dots
\end{itemize}

Clearly, there are other important groups e.g. in {\bf Vienna} (e.g. J. Schoissengeier, J. Schwermer, C. Baxa, L. Summerer, M. Drmota, G. Dorfer, E. Kr\"atzel, W.G. Nowak, M. K\"uhleitner, K. Scheicher, P. Surer, \dots) and {\bf Linz} (G. Larcher, F. Pillichshammer, H. Niederreiter, A. Winterhof, \dots); we will not mention their important contributions in number theory in this paper (unless they were part of the above groups in the past, as K. Scheicher, who did his PhD in Graz with R. Tichy, or P. Surer, who did his PhD in Leoben with J. Thuswaldner).\footnote{For example we will not further mention the important recent joint work of W.M. Schmidt and L. Summerer (cf. \cite{MR2557854,MR3016519,MR3284111}), where they achieved a considerable progress on old questions of Jarnik using their parametric geometry of numbers, even if these results clearly belong to Diophantine number theory.}

Because of his obvious special importance to the Austrian number theory community, we also mention W.M. Schmidt (University of Colorado, Boulder/US). Much of our research is inspired or motivated by his work.\\

Now we turn to the Hungarian research group in which the second author is involved. It is the Number Theory and Cryptography School of Debrecen, whose founder and leader is K. Gy\H{o}ry. Its members include (given in order of receiving their PhD degrees / starting their PhD studies):
\begin{itemize}
\item {\bf Debrecen}: \dag B. Kov\'acs, A. Peth\H{o}, S. Turj\'anyi, \dag Z. Papp,\break \dag B. Brindza, I. Ga\'al, \'A. Pint\'er, L. Hajdu, A. B\'erczes, T. Herendi,\break Sz. Tengely, I. Pink, G. Nyul, A. Huszti, A. Bazs\'o, J. Foll\'ath,\break Gy. P\'eter, Zs. R\'abai, N. Varga, J. Ferenczik, B. Rauf, E. Gyimesi,\break T. Szab\'o, M. Szikszai, Cs. Bert\'ok, G. R\'acz
\item {\bf Eger}: \dag P. Kiss, F. M\'aty\'as, K. Liptai, T. Szak\'acs
\item {\bf Marosv\'as\'arhely (Tirgu Mures)}: Gy. M\'arton
\item {\bf Miskolc}: Cs. Rakaczki, P. Olajos
\item {\bf Ny\'\i regyh\'aza}: J. K\"odm\"on, Z. Csajb\'ok
\item {\bf Sopron}: L. Szalay, K. Gueth
\item {\bf Budapest}: T. Kov\'acs
\end{itemize}
The group is ``vivid'' and not at all closed. Members move in and out of the research group (e.g. to focus more on cryptography) at certain times.

Certainly, (beside rich international contacts) the group has strong connections with many other excellent Hungarian researchers and research groups in number theory (e.g. from the R\'enyi Institute, E\"otv\"os University (ELTE), TU Budapest, Institute for Computer Science and Control (SZTAKI), etc.).\\

Having localized now the participants of our collaboration, we shall give some interesting data of our groups in the next section.

\section{Some of the most important joint research topics and results}

The main joint research topics practically cover all parts of Diophantine number theory, and several topics from representation theory, theory of (recurrence) sequences, diophantine approximation, algebraic number theory and exponential sums. We tried to extract some interesting data from Mathematical Reviews (MR). This turned out to be not that easy, since the data in MR sometimes lacked of completeness, e.g. sometimes there are different affiliations attached to one and the same author or sometimes there is no affiliation attached at all. (We mention P. Kiss, A. Peth\H{o} and also the first author as examples for such assignment problems.) Thus the data below has to be taken with some care and should be seen as verified lower bounds. We have included the list of references, taken from MR, to which we refer to; see \cite{MR2075644} up to \cite{MR3327847}.\\

The list of {\bf primary} MR classification numbers, with joint publications are the following: \begin{quote}05C25, 11A63, 11A67, 11B13, 11B30, 11B37, 11B75, 11D09, 11D25, 11D41, 11D45, 11D57, 11D59, 11D61, 11D72, 11J04, 11J25, 11K06, 11K16, 11K36, 11K60, 11K70, 11R06, 11R58, 11T23.\end{quote} Here some more details on those 2010 Mathematical Subject Classification numbers, which appear most often, ordered by the number of appearances:

\begin{itemize}
\item[11A63] Elementary number theory. Radix representation; digital problems - 10 times
\item[11K16] Probabilistic theory: distribution modulo $1$; metric theory of algorithms. Normal numbers, radix expansions, Pisot numbers, Salem numbers, good lattice points, etc. - 6 times
\item[11B37] Sequences and sets. Recurrences - 5 times
\item[11D41] Diophantine equations. Higher degree equations; Fermat's equation - 5 times
\item[11D59] Diophantine equations. Thue-Mahler equations - 5 times
\item[11D45] Diophantine equations. Counting solutions of Diophantine equations - 4 times
\item[11D09] Diophantine equations. Quadratic and higher degree equations - 3 times
\end{itemize}

According to MR, the two groups have at least 64 joint publications so far in the above mentioned topics. The recorders are A. Peth\H{o} with 36 and R. Tichy with 24 joint papers. They have 15 joint papers as coauthors. We also mention J. Thuswaldner with at least 11 joint papers, V. Ziegler with at least 10 joint papers, the first author with 7 joint papers with Hungarian colleagues as well as P. Kiss with at least 5 joint papers, \'A. Pint\'er and L. Szalay both with at least 4 joint papers and A. B\'erczes with at least 3 joint papers with Austrian colleagues from our joint collaboration.\\

To the best of our knowledge, the first joint publication is (we will state one of the results from that paper in Subsection \ref{lrs} below): \begin{quote}P. Kiss, R. Tichy, Distribution of the ratios of the terms of a second order linear recurrence. \emph{Nederl. Akad. Wetensch. Indag. Math.} {\bf 48} (1986), no. 1, 79--86.\medskip\end{quote}

Below we shall give some more details on joint research topics. We stress that we will concentrate on {\bf joint} research. Certainly, there are several other topics where the members of the research groups obtained important (in many cases outstanding) results. The topics which will be highlighted here are:
\begin{itemize}
\item[3.1] Parametric families of Thue equations
\item[3.2] Canonical number systems and shift radix systems
\item[3.3] Separated-variables Diophantine equations
\item[3.4] Linear recurrence sequences
\end{itemize}
The details will be given in the next four subsections. In each of the subsections we shall first introduce some important relevant notions and give main references before we highlight some results from our collaboration from the last 30 years.

\subsection{Parametric families of Thue equations}

Let $f\in\mathbb{Z}[X,Y]$ be an irreducible form of degree at least $3$ and $m\in\mathbb{Z}\backslash\{0\}$. The equation \[f(x,y)=m\]
in $x,y\in\mathbb{Z}$, is called a \emph{Thue equation}. Some known results in chronological order:
\begin{itemize}
\item {Thue \cite{MR1580770}}: such a (Thue) equation has only finitely many solutions,
\item {Baker \cite{MR0228424}, Bugeaud and Gy\H{o}ry \cite{MR1373714}}: the solutions can be calculated effectively,
\item {Baker and Davenport \cite{MR1580770}, Peth\H{o} and Schulenberg \cite{MR0934900}, Tzanakis and de Weger \cite{MR0987566,MR1094961,MR1151871,MR1189890,MR1255696}, Bilu and Hanrot \cite{MR1412969}}: the solutions can be calculated efficiently\footnote{Because of the huge bounds that typically arise by applying Baker's method of linear forms in logarithms of algebraic numbers, calculating all solutions to a conrete Thue equation was - and still is - a non-trivial task; several people worked on this problem.},
\item {Gill \cite{MR1502935}, Schmidt \cite{MR0508466}, Mason \cite{MR0635873}, Dvornicich and Zannier \cite{MR1309649}}: function field analogues of the above results,
\item {Gy\H{o}ry \cite{MR0716553}, Gy\H{o}ry, B\'erczes and Evertse \cite{MR3194058}}: results for finitely generated domains.
\end{itemize}

Since it is known that a single Thue equation has only finitely many solutions and these solutions can, at least in principle, be calculated, the research focused on solving families of Thue equations \[f_a(x,y)=m\] in $x,y\in\mathbb{Z}$, with $f_T\in\mathbb{Z}[T][X,Y],a\in\mathbb{Z}$. Results are available for (for a survey see \cite{heuberger}):
\begin{itemize}
\item families of fixed degree, e.g. by {Thomas\footnote{This paper can be seen as starting point of the investigations on parametrized Thue equations.} \cite{MR1042497}, \dots
\footnote{We are not giving more references here but refer to \cite{heuberger} instead (even if that list is certainly not up to date anymore, e.g. \cite{MR2271282} and \cite{MR2492925} are two examples of papers on this topic which appeared afterwards). References with contributions from individual members from our research groups are e.g. \cite{MR1770464,MR2093887,MR1397316,MR1094956,MR1316596,MR1359142,MR1487624,MR1628852,MR1635234,MR1822519,MR1635238,MR1659855,MR1796800,MR2093887}.},
Chen, Ga\'al, Heuber\-ger, Jadrijevi\'c,  Lee, Lettl, Levesque, Mignotte, Peth\H{o}, Tichy, Togb\'e, Roth, Tzanakis, Voutier, Wakabayashi, Ziegler},
\item families of relative equations, e.g. by {Heuberger, Peth\H{o} and Tichy \cite{MR1937468,MR2373935}, Heuberger \cite{MR2251815}, Ziegler \cite{MR2194942,MR2242388,MR2684373}, Jadrijevi\'c and Ziegler \cite{MR2281865}},
\item families of arbitrary degree, e.g. by {Halter-Koch, Lettl, Peth\H{o} and Tichy \cite{MR1721811}, Heuberger \cite{MR1688196,MR1829779,MR1844070}, Heuberger and Tichy \cite{MR1726188}},
\item function field analogues, e.g. by {Fuchs and Ziegler \cite{MR2199120,MR2204261}, Lettl \cite{MR2114617,MR2139595}, Ziegler \cite{MR2332067,MR2646467}, Fuchs and Jadrijevi\'c \cite{MR2340972}, Ga\'al and Pohst \cite{MR2228949,MR2229380}}, Fuchs, Jurasi\'c and Paulin \cite{arXiv2}.
\end{itemize}
There are already many names from our research groups mentioned in this list. We give four explicit examples of joint results:
\begin{itemize}
\item {Halter-Koch, Lettl, Peth\H{o} and Tichy (\cite{MR1721811}; 1996)}: Let $n\geq 3$, $a_1=0,a_2,\ldots,a_{n-1}$ be distinct integers and $a_n=a$ an integral parameter. Let $\alpha=\alpha(a)$ be a zero of $P(x)=\prod_{i=1}^n(x-a_i)-d$ with $d=\pm 1$ and suppose that the index $I$ of $\langle \alpha-a_1,\ldots,\alpha-a_{n-1}\rangle$ in the group of units of $\mathbb{Z}[\alpha]$ is bounded by a constant $J=J(a_1,\ldots,a_{n-1},n)$ for every $a$ from some subset $\Omega\subseteq\mathbb{Z}$. Assume further that the Lang-Waldschmidt conjecture is true. Then for all but finitely many values $a\in\Omega$ the equation $\prod_{i=1}^n(x-a_iy)-dy^n=\pm 1$ has only solutions $(x,y)\in\mathbb{Z}^2$ with $\vert y\vert\leq 1$, except when $n=3$ and $\vert a_2\vert=1$, or when $n=4$ and $(a_2,a_3)\in\{(1,-1),(\pm 1,\pm 2)\}$ in which case there is exactly one more solution for every values of $a$.
\item {Heuberger, Peth\H{o} and Tichy (\cite{MR1822519}; 1998)}: The family of Thue equations given by $x(x-y)(x-ay)(x-(a+1)y)-y^4=\pm 1$ has only the trivial solutions $(x,y)\in\{(\pm 1,0),(0,\pm 1),$ $(\pm 1,\pm 1),(\pm a,\pm 1),(\pm (a+1),\pm 1)\}$ in $\mathbb{Z}^2$.
\item {Ga\'al and Lettl (\cite{MR1659855,MR1796800}; 2000)}: The family of Thue equations given by $x^5+(a-1)^2x^4y-(2a^3+4a+4)x^3y^2+(a^4+a^3+2a^2+4a-3)x^2y^3+(a^3+a^2+5a+3)xy^4+y^5=\pm 1$ has only trivial solutions $\pm(x,y)\in\{(1,0),(0,1)\}$ for $a\in\mathbb{Z}\backslash\{-1,0\}$ and $(x,y)\in\{(1,0),(0,1),(\pm 1,1),(-2,1)$ for $a=-1,0$, resp., in $\mathbb{Z}^2$.
\item {Fuchs, Jurasi\'c and Paulin (\cite{arXiv2})}: The family of Thue equations given by $x(x-y)(x+y)(x-ay)+y^4=\xi$ with $\xi\in\mathbb{C}^\times$ and $a\in\mathbb{C}[T]\backslash\mathbb{C}$ has only the trivial solutions $(x,y)\in\{(\zeta,0),(0,\zeta),(\zeta,\zeta),(-\zeta,\zeta),$ $(\zeta a,\zeta),(-\zeta,\zeta a);\zeta\in\mathbb{C}^\times,\zeta^4=\xi\}$ in $\mathbb{C}[T]^2$.
\end{itemize}

More results can be found e.g. in \cite{MR1937468,MR2373935,MR1359142,MR1487624,MR1628852}.

\subsection{Canonical number systems and shift radix systems}

Number systems and their generalizations have been studied for a very long time.

As it is well-known, any positive integer $n$ can be uniquely written as
\[
n=\sum\limits_{i=0}^k d_ib^i,
\]
where $b\geq 2$ is a fixed integer (the base of the number system), and $d_i\in\{0,1,\dots,b-1\}$ $(d_k\neq 0)$ are the digits.

{Gr\"unwald \cite{Grunwald1885}} showed that this can be extended to the ring of integers, using negative base $b$.

Extensions to rings of integers of algebraic number fields were worked out by Knuth \cite{MR0127508}, Penney \cite{Penney1965}, K\'atai and Szab\'o \cite{MR0389759}, K\'atai and Kov\'acs \cite{MR0576942, MR0616887}, Kov\'acs \cite{MR0619892}, Gilbert \cite{MR0632342}, Kov\'acs and Peth\H{o} \cite{MR1152592}, K\'atai and K\"ornyei \cite{MR1189110}, and many others. The main directions of research concern the description of rings having a canonical number system, and the description of the canonical number systems if they exist. The results are partly based upon deep theorems of Gy\H{o}ry \cite{MR0784078} yielding effective upper bounds for generators of power integral bases in number fields.

{Peth\H{o} (\cite{MR1151853}; 1991)} extended the notion to residue class rings of polynomials (called CNS polynomials). If $p(x)$ is irreducible with a root $b$, then $\mathbb{Z}[b]$ is isomorphic to $\mathbb{Z}[x]/p(x)\mathbb{Z}[x]$. So the two notions are very much interrelated.

Such CNS polynomials were then described in various cases by {Kov\'acs \cite{MR0619892}, K\'atai and Kov\'acs \cite{MR0616887}, K\"ormendi \cite{MR0882046}, Brunotte \cite{MR1876451}, Akiyama and Peth\H{o} \cite{MR1871104}, Scheicher and Thuswaldner \cite{MR2069090} and others; see also the survey paper \cite{MR2075644} and the references there.

Later, Akiyama, Borb\'ely, Peth\H{o} and Thuswaldner (\cite{MR2162561}; 2005)} made the following generalization:

Let ${\mathbf r}=(r_1,\dots,r_d)\in {\mathbb R}^d$ and $\tau_{\mathbf r}:\ {\mathbb Z}^d\to {\mathbb Z}^d$ be defined by $\tau_{\mathbf r}(a_1,\dots,a_d)= (a_2,\dots,a_d,-\lfloor r_1a_1+\dots+r_da_d\rfloor).$ If for all ${\mathbf a}$ we have $\tau_{\mathbf r}^k({\mathbf a})={\mathbf 0}$ for some $k$, then $\tau_{\mathbf r}$ is a SRS (short for shift radix system). Put $D_d^{0}=\{{\mathbf r}\ :\ \tau_{\mathbf r}\ \text{is SRS}\}.$\footnote{The elements of $D_d^{0}$ are called ``SRS with finiteness property'' in the recent papers, which indeed seems to be a more adequate name.}

The importance of shift radix systems is shown by the fact that $p(x)=x^d+p_{d-1}x^{d-1}+p_1x+p_0$ $(p_i\in{\mathbb Z})$ is CNS if and only if $(1/p_0,p_{d-1}/p_0,\dots,p_1/p_0)\in D_d^{0}$; see \cite{MR2162561, MR2735753}.

In the past few years, extensive research has been done and interesting results concerning SRS-s have been obtained, by Akiyama, Borb\'ely, Brunotte, Peth\H{o} and Thuswaldner \cite{MR2162561}, Akiyama, Brunotte, Peth\H{o} and Thuswaldner \cite{MR2216302, MR2441944, MR2513054}, Akiyama, Brunotte, Peth\H{o} and Steiner \cite{MR2224891}, Akiyama, Brunotte, Peth\H{o}, Steiner and Thuswaldner \cite{abpst}, Berth\'e, Siegel, Steiner, Surer and Thuswaldner \cite{MR2735753}, Huszti, Scheicher, Surer and Thuswaldner \cite{MR2342432},
Berth\'e, Siegel, Steiner, Surer and Thuswaldner \cite{MR2735753}, Brunotte, Kirschenhofer and Thuswaldner \cite{MR2907970}, Madritsch and Peth\H{o} \cite{MR3129005}, Weitzer \cite{MR3367922} and Peth\H{o}, Varga and Weitzer \cite{pvw}. See also the survey papers \cite{MR2290774, MR3330559}.

To show a particular, interesting example, extending the problem to the complex field, {Brunotte, Kirschenhofer, Thuswaldner (\cite{MR2907970}; 2011)} for the complex analogue of $D_1^{0}$ obtained Peth\H{o}'s loudspeaker\footnote{Notice that $D_1^{0}$ is trivial, the interval $(-1,1)$, but its complex analogue is the loudspeaker.}$^,$\footnote{The picture is taken from \cite{MR2907970} with the permission of its authors.}:
\begin{figure}[htp!]
  \begin{center}
    {\includegraphics[scale=0.25]{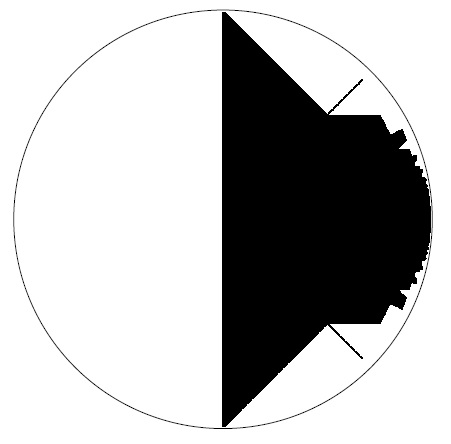}}
  \end{center}
\end{figure}

\subsection{Separated-variables Diophantine equations}

Let $f,g\in\mathbb{Z}[T]$. A Diophantine equation of the form \[f(x)-g(y)=0\] in integers $x,y$ is called \emph{of separated-variables type}.

We start by mentioning the most important general result(s) on Diophantine equations of separated-variables type\footnote{Observe that we are not giving a full historical account; for this one should also mention (at least) results of Davenport, Lewis and Schinzel (cf. \cite{MR0137703} and the monograph \cite{MR1770638}) and the important contributions due to Fried (cf. the survey article \cite{MR2873803}).}:
\begin{itemize}
\item The {Bilu-Tichy} theorem \cite{MR1793164} says that such an equation has infinitely many solutions (in rationals with bounded denominator) if and only if $f=\varphi\circ f_1\circ \lambda,g=\varphi\circ g_1\circ\mu$ with $\varphi,\lambda,\mu\in\mathbb{Q}[T]$, $\deg(\lambda)=\deg(\mu)=1$ and $(f_1,g_1)$ is of special (precisely described) form. Thus to prove that such a Diophantine equations has only finitely many solutions, one just has to verify that no such decomposition exists.
\item {Bilu, Fuchs, Luca and Pint\'er (\cite{MR3034377}; 2013)}: In the above statement almost all solutions satisfy $f_1(\lambda(x))=g_1(\mu(y))$; the possible exceptions can be completely described.
\end{itemize}

The Bilu-Tichy theorem has many applications. We mention a few of them:
\begin{itemize}
\item {Bilu, Brindza, Kirschenhofer, Pint\'er and Tichy (\cite{MR1898434}; 2002)}: $x(x+1)\cdots(x+m-1)=1^n+2^n+\cdots+(x-1)^n$ has only finitely many solutions for $m\geq 2,n\geq 1, (m,n)\neq (2,1).$
\item {Stoll and Tichy (\cite{MR2027780,MR2035894,MR2136185,MR2372822}; 2003-2008)}: Finiteness results for Diophantine equations involving various orthogonal (classical continuous and modified orthogonal) polynomials, binomial coefficients, general Meixner and Krawtchouk polynomials.
\item {Bilu, Fuchs, Pint\'er and Luca (\cite{MR3034377}; 2013)}: The Diophantine equations $S_{x-a}^x=S_{y-b}^y,s_{x-a}^x=s_{y-b}^y,$ $S_{x-a}^x=s_{y-b}^y,S_{x-a}^x={y\choose b}, s_{x-a}^x={y\choose b}$ have only finitely many solutions $(x,y)\in\mathbb{N}^2$, where $S$ denotes the Stirling numbers of the second and $s$ the (unsigned) Stirling numbers of the first kind.
\item {Kreso and Rakaczki (\cite{MR3145451}; 2013)}: Finiteness results for Diophantine equations with Euler polynomials.
\end{itemize}

Yet another example is the ``Debrecen problem'': Let $n\geq 1,r\geq 1$ and \[f_n(x)=\#\{(x_1,\ldots,x_n)\in\mathbb{Z}^n\ :\ \vert x_1\vert+\cdots+\vert x_n\vert\leq x\}.\]
It was proved:
\begin{itemize}
\item {Hajdu (\cite{MR1485228,MR1604054}; 1997-1998)}:  The equations $f_2(x)=f_3(y)$ and $f_2(x)=f_4(y),f_2(x)=f_6(y),f_3(x)=f_4(y),$ $f_4(x)=f_6(y)$ have only finitely many solutions.
\item {Kirschenhofer, Peth\H{o} and Tichy (\cite{MR1702180}; 1999)}: $f_n(x)=f_m(y)$ has only finitely many solutions for $n=2,m\geq 3$ or $n=4,m\geq 3$ or $n,m\geq 2$ not congruent modulo $2$ or $2\leq n<m\leq 103$.
\item {Bilu, Stoll and Tichy (\cite{MR1793164}; 2000)}: $f_n(x)=f_m(y)$ has only finitely many solutions for any pair $(n,m)$ of distinct positive integers with $n,m\geq 2$.
\end{itemize}

This problem originates from the investigations of the second author and was brought to Graz by A. Peth\H{o}. There it was named the ``Debrecen problem'' by R. Tichy. This is a very nice example of our collaboration.

\subsection{Linear recurrence sequences}\label{lrs}

A sequence $(G_n)$ (in some ring, to be specified later) is called a \emph{linear recurrence} if the terms satisfy an equation of the form $$G_{n+d}=a_1G_{n+d-1}+\cdots+a_dG_n$$ for all $n\geq d$ for given $a_1,\dots,a_d$ and $G_0,\ldots,G_{d-1}$. If $\alpha_1,\ldots,\alpha_t$ denote the distinct characteristic roots of the sequence (i.e. of the polynomial $x^n-a_1x^{n-1}-\dots-a_d$), then we have the formula $G_n=P_1(n)\alpha_1^n+\cdots+P_t(n)\alpha_t^n$, where $P_i$ are polynomials with coefficients in the splitting field of the characteristic equation of degree less than the multiplicity of the corresponding characteristic root.

A recurrence is called \emph{non-degenerate} if no quotient $\alpha_i/\alpha_j$ is a root of unity and it is called \emph{non-unitary} if no characteristic root $\alpha_i$ is a root of unity.

Many arithmetical properties have been studied for linear recurrences including e.g. digital and Diophantine aspects.\footnote{Instead of trying to review the huge related literature, we only refer to the books \cite{MR0891406, MR1990179}.}

The first joint result of the research groups by {Kiss and Tichy (\cite{MR0834322}; 1986)} was the following:
Let $(G_n)$ be a real second order linear recurrence satisfying $G_{n+2}=AG_{n+1}+BG_n$ with discriminant $D=A^2+4B<0$ and assume that $\pi^{-1}\arctan(\sqrt{-D}/A)$ is irrational. Then the fractional parts of $G_{n+1}/G_n$ is everywhere dense in $[0,1]$ but is not uniformly distributed.

We mention also the following results, which all belong to joint papers from our collaboration:
\begin{itemize}
\item {Kiss and Tichy (\cite{MR0906995,MR1011853,MR1011869}; 1987, 1989), Grabner, Kiss and Tichy (\cite{MR1353931}; 1995)}: Results concerning the distribution and approximation properties of ratios of consecutive terms of binary linear recurrence sequences.
\item {Peth\H{o} and Tichy (\cite{MR1034204,MR1208858}; 1989, 1993)}: Results for the sum-of-digits function for representations of integers with respect to a linear recurrence sequence.
\item {Fuchs, Peth\H{o} and Tichy (\cite{MR1942618,MR1990766,MR2487694}; 2002-2008)}:
Let $(G_n)$ be a non-degenerate recurrence sequence of polynomials and let $P$ be a polynomial. Under (fairly general) conditions\footnote{In \cite{MR2111014} the conditions were completely clarified by Zannier and a sharp bound on the number of solutions was given.} it follows that $G_n(x)=G_m(P(x))$ has only finitely many solutions $(n,m)\in\mathbb{Z}^2$ with $n,m\geq 0,n\neq m$.
\item {Fuchs and Peth\H{o} (\cite{MR2212123}; 2005)}: Let $(G_n)$ be a non-degenerate recurrence which is ``truely'' defined over a function field with characteristic zero constant field. Then there is an effectively computable constant $C$ such that for all $n\in\mathbb{N}$ with $G_n=0$ we have $n\leq C$. In the same paper a similar result was shown to hold for $G_n=H_m$ and therefore also for $G_n(x)=G_m(P(x))$.
\item {Fuchs, Luca and Szalay (\cite{MR2483637}; 2008)}: Let $(G_n)$ be a non-degenerate binary recurrence with positive discriminant such that for infinitely many sextuples of non-negative integers $(a,b,c;x,y,z)$ with $1\leq a<b<c$ we have $ab+1=G_x,ac+1=G_y,bc+1=G_z$. Set $G_n=\gamma\alpha^n+\delta\beta^n$. Then $\beta,\delta\in\{\pm 1\},\alpha,\gamma\in\mathbb{Z}$. Indeed, for recurrences of this form (e.g. $G_n=2^n+1$) infinitely many such sextuples exist!
\item {Pint\'er and Ziegler (\cite{MR2922338}; 2012)}: If a non-unitary, non-degenerate binary recurrence contains infinitely many three-term arithmetic progressions, its companion polynomal divides $(X^a-2X^b+1)/(X^d-1)$ or $(X^a+X^b-2)/(X^d-1)$ with some $a>b>0$, $d=\gcd(a,b)$. This implies that the only increasing, non-unitary, non-degenerate recurrence with initial terms $0,1$ and $\alpha\neq\beta$ containing infinitely many three term arithmetic progressions is the Fibonacci sequence.
\item {B\'erczes and Ziegler (\cite{MR3064238}; 2013)}: Let $G_n$ be a Lucas sequence (i.e. a binary recurrence with $G_0=0$, $G_1=1$). Then $G_n$ does not contain non-trivial geometric progressions, except for the one characterized by $(G_1,G_2,G_4)= (G_3,G_2,G_4)=(1,-2,4)$.
\end{itemize}

More results can be found for example in \cite{MR3081229,MR3167927,MR3327847,arXiv1}. Other lines of related research which we will not go into can e.g. be found in \cite{MR2232201,MR2419169}.


\section{A closer look - on a conjecture by R\'enyi and Erd\H{o}s}

In view of the mentioned applications to Diophantine equations of separated-variables type, it is interesting to ask for results on the functional decomposition of polynomials and rational functions. For larger classes of polynomials this is a highly non-trivial problem. Below we shall restrict to fewnomials (also known as lacunary resp. sparse polynomials). Let $g(x)\in\C[x]$ be given. We denote by $\N(g)$ its number of terms, i.e. we have \[g(x)=a_1x^{n_1}+\cdots +a_{\N(g)}x^{n_{\N(g)}},\] wehre $a_j\neq 0$ for $j=1,\ldots,\N(g).$
Now given $\ell\in\mathbb{N}$. All polynomials in the set \[\{g(x)\in\C[x]\ :\ \N(g)\leq \ell\}\] are called fewnomials/lacunary resp. sparse polynomials (with respect to $\ell$).\\

An attractive conjecture was formulated independently by R\'enyi and Erd\H{o}s (we mention the papers \cite{MR0022268} of R\'enyi from 1947 and \cite{MR0027779} of Erd\H{o}s from 1949): Is there a $B_0=B_0(\ell)$ such that if $\N(g^2)\leq \ell$, then $\N(g)\leq B_0(\ell)$? In other words: If $g^2=\sum_{i=1}^\ell a_ix^{n_i}$, does it follow that $\N(g)\leq B_0$ for some $B_0=B_0(\ell)$?

This was indeed proved by Schinzel in 1987 (see \cite{MR0913764}) in more general form and with an explicit value for the constant\footnote{Schinzel's result even holds for polynomials in a field of positive characteristic. The results below do not have this feature; for the purpose of this presentation we always stick to the case of complex polynomials.}$^,$\footnote{The result in \cite{MR0913764} says that $\ell\geq d+1+(\log 2)^{-1}\log(1+\log(\N(g)-1)/(d\log(4d)-\log d))$; we give a slightly simplified version here.}:
If $\N(g^d)\leq \ell$, then $\N(g)\leq 1+$ $(4d)^{2^{\ell-d-1}d}$.

Schinzel conjectured that even more might be true: If $\deg(f)\leq d$ and $\N(f(g(x)))\leq \ell$, then $\N(g)\leq B_1(d,\ell)$ for some $B_1=B_1(d,\ell)$.

This conjecture was proved by Zannier (cf. \cite{MR2430978}) in 2008. Indeed one can take: \[B_1=(4\ell)^{(2\ell)^{(3\ell)^{\ell+1}}}\] (as was recently shown by Karolus \cite{karolus}). We do not claim that the bound is sharp, however, we do not have good lower bounds. In fact for the additional uniformity in the last result one needs the following, which Zannier proved in 2007 (see \cite{MR2289981} and also \cite{MR2557855}):
If $\N(f(g(x)))\leq \ell$, then either $g(x)=ax^n+b$ or $\deg(f)\leq$ $\ell(\ell-1)$. Thus one can take $B_1=B_1(\ell)$.\\

We now turn to rational functions; this is quite natural if one thinks of Diophantine problems involving linear recurrences, where typically families of solutions come from identities for Laurent polynomials.  Let $g(x)\in\C(x)$ be given. We define \[\N^\#(g)=\min\{\N(p)+\N(q)\ :\ g=p/q\}.\]
Note that $p,q$ may not be coprime. Thus $\N^\#(1+\cdots+x^{n-1})$ $=\N^\#((x^n-1)/(x-1))=4$. Therefore $\N^\#$ does not coincide with $\N$ on polynomials. This shows that indeed the fact that we do not necessarily consider coprime representations for $g$ makes the investigation of the number of terms much more involved.

The first author and Zannier proved in 2012 (cf. \cite{MR2862037}): If $\N^\#(f(g(x)))\leq \ell$, then either $g(x)=\lambda(ax^n+bx^{-n})$ ($\lambda\in$ PGL$_2(\mathbb{C})$, $a,b\in\mathbb{C},n\in\mathbb{N}$) or $\deg (f)\leq 2016\cdot 5^\ell$. This is the analogue of Zannier's result from \cite{MR2289981} and shows that rational operations ``tend to destroy'' the lacunarity.\\

The following question rises: What about $\N^\#(g)$ and thus Schinzel's conjecture? A recent result is the following:

\begin{theorem}[Fuchs, Mantova and Zannier \cite{arXiv3}] If $\deg(f)\leq d$ and $\N^\#(f(g(x))\leq \ell$, then $\N^\#(g)\leq B_2=B_2(d,\ell)$.
\end{theorem}

A combination with the result from \cite{MR2862037} gives $\N^\#(g)\leq B_3(\ell)$. This can be reformulated in the following way:
Given $\ell\in\mathbb{N}$. There exist the effectively computable data $J,t\in\mathbb{N}$ and for all $1\leq j\leq J$: $\mathcal{V}_j/\mathbb{Q}$ and $F_j,H_j\in\mathbb{Q}[\mathcal{V}_j](x_1,\ldots,x_t),G_j\in$ $\mathbb{Q}[\mathcal{V}_j](x)$ with $F_j=G_j\circ H_j$ and $F_j$ has $\ell$ terms such that the following holds: \begin{quote}\textnormal{If $f(x)=g(h(x)),g,h\in \mathbb{C}(x)$ has $\ell$ terms, $h(x)\neq$ $\lambda(ax^n+bx^{-n})$, then there is $j$ and $P\in\mathcal{V}_j(\CC),$ $(k_1,\ldots,k_t)\in\mathbb{Z}^t$ with $g(x)=G_j(P,x), h(x)=H_j(P,x^{k_1},\ldots,x^{k_t}), f(x)=F_j(P,x^{k_1},\ldots,x^{k_t}).$}\end{quote}

But we can do much more! We have shown that rational roots of lacunary equations are lacunary, in more details the following result is true.

\begin{theorem}[Fuchs, Monatova and Zannier \cite{arXiv3}] Let $f(x,y)\in\C[x,y]$ with $\N_x(f)\leq \ell$ and $\deg_y(f)\leq d$. \begin{itemize} \item[a)] If $f(x,g(x))=0$ with $g(x)\in\C(x)$, then $\N^\#(g)\leq B_4=B_4(d,\ell).$
\item[b)] If $f$ is monic, then $g(x)\in\C[x]$ and $\N(g)\leq B_5=B_5(d,\ell).$
\end{itemize}
\end{theorem}

The result can be reformulated in the following obvious equivalent way, which is more suitable for the proof and also for applications:

\begin{theorem}[Fuchs, Mantova and Zannier \cite{arXiv3}] Let $f(x_1,\ldots,x_\ell,y)\in\C[x_1,\ldots,$ $x_\ell,y]$ and $f$ have degree at most $d$ in each variable. \begin{itemize} \item[a)] If $f(x^{n_1},\ldots,x^{n_\ell},g(x))=0$ with $g(x)\in\C(x)$, then $\N^\#(g)\leq B_6=B_6(d,\ell).$ \item[b)] If $f$ is monic in the last variable, then $g(x)\in\C[x]$ and $\N(g)\leq B_7=B_7(d,\ell).$\end{itemize}
\end{theorem}

The constants $B_2$ to $B_7$ can all be calculated from the proof. However, they will be quite poor (certainly iterated exponential) so that no explicit values were worked out.

We shall not go into mentioning applications of these results (which can be found in \cite{arXiv3}) or elaborate too much on the proof. Instead, we just mention a few ingredients of it. The proof is done by an intricate triple induction argument; the main induction is done on $\ell$. Without loss of generality we may assume that $f$ is monic in $y$. Moreover, the proof uses:
\begin{itemize}
\item Newton polygons and Puiseux series (in a more refined form) - we have found a new way of making the expansion which gives additional information,
\item Puiseux series of other base fields than the rational function field - this was not available in the literature and might be of interest in its own,
\item Pseudo-Puiseux expansion, which we have developed to control the indices in the Puiseux expansions,
\item successive minima for the exponents which gives a reparametrization procedure,
\item Diophantine approximation in function fields.
\end{itemize}

To finish with this part, we mention that similar results can be given for other types of lacunarity. One can consider families of polynomials/rational functions which are inspired by varying the concept of lacunarity, e.g.:
\begin{itemize}
\item fixing the number of zeros and poles; this was considered by {Fuchs and Peth\H{o} (\cite{MR2729068}; 2011)}\footnote{See also \cite{MR3220123} for some further information in this direction.},
\item taking polynomials satisfying a common linear recurrence; special cases of this were considered by {Fuchs, Peth\H{o} and Tichy (\cite{MR1942618,MR2040240,MR1990766,MR2487694}; 2002-2008)},
\item clearly, there are many more possibilities.
\end{itemize}


\section{A closer look - the unit sum number problem}

An integral domain $R$ is called $k$-good, if every element of $R$ can be written as a sum of $k$ units. The unit sum number $u(R)$ of $R$ is:
\begin{itemize}
\item the minimal integer $k$ such that $R$ is $k$-good, if such an integer exists,
\item $\omega$, if $R$ is not $k$-good for any $k$ but every element of $R$ is a sum of units,
\item $\infty$ otherwise, i.e. if $R$ is not (additively) generated by its units.
\end{itemize}

Some simple examples:
\begin{itemize}
\item $u({\mathbb Q})=2$,
\item $u({\mathbb Z})=\omega$,
\item $u({\mathbb Q}[x])=\infty$.
\end{itemize}

There are several results concerning different kinds of domains $R$, for example:
\begin{itemize}
\item endomorphism rings of modules and vector spaces: {Zelinsky (\cite{MR0062728}; 1954)}, {Goldsmith, Pabst and Scott (\cite{MR1645560}; 1998)},
\item matrix rings: {Henriksen (\cite{MR0349745}; 1974), Levy (\cite{MR0294367}; 1972)}, {V\'amos and Wiegand (\cite{MR2806699}; 2011)},
\item function fields: {Frei (\cite{MR2827171,MR2818674}; 2011)},
\item rings of integers of quadratic number fields and cyclotomic fields: {Ashrafi and V\'amos (\cite{MR2124574}; 2005)}.
\end{itemize}

Now we shall focus exclusively on rings of integers of algebraic number fields. Let $K$ be an algebraic number field, with unit group $U$. For given $k\geq 2$, consider the so-called unit equation
\[
u_1+\dots+u_k=1
\]
in unknowns $u_1,\dots,u_k\in U$. A solution to this equation is {\sl non-degenerate}, if no subsum on the left hand side vanishes. There are many deep results concerning unit equations; most importantly, such an equation admits only finitely many non-degenerate solutions.\footnote{Once again, instead of giving a survey of the vast literature, we only refer the interested reader to the excellent survey papers \cite{MR0971998,MR1211001} and books \cite{MR0891406,egy2015}.} This is ultimately based upon the famous subspace theorem of {Schmidt} (cf. e.g. \cite{MR0568710,MR1176315}).

Unit equations appear in a huge variety of applications, e.g. to Diophantine equations, irreducibility of polynomials, to arithmetic graphs, properties of digital sequences like discrepancy of lacunary sequences, etc.\footnote{See again \cite{MR0971998,MR1211001,MR0891406,egy2015}.}

Let $K$ be a number field, with unit group $U$. Write $r$ for the rank of $U$. Let $k$ be a positive integer, and set
\[
H_k=\left\{u_1+\dots+u_k\
:\ u_1,\hdots,u_k\in U\right\}.
\]

\begin{theorem}[Jarden, Narkiewicz \cite{MR2309537}, Hajdu \cite{MR2337316}, independently] The length of any non-constant arithmetic progression in $H_k$ is at most $C(r)$.
\end{theorem}

The main tools in the background of the proof are:
\begin{itemize}
\item a result of {Evertse, Schlickewei and Schmidt \cite{MR1923966}} on the non-degenerate solutions of unit equations,
\item a classical result of {Van der Waerden \cite{MR0175875}} on monochromatic arithmetic progressions (to handle vanishing subsums).
\end{itemize}

The result in fact holds more generally on replacing $U$ by a multiplicative subgroup of the multiplicative group of $K$ of, say, rank $r$.

\begin{corollary}[Jarden and Narkiewicz \cite{MR2309537}] The ring $R$ of integers of any algebraic number field $K$ is not $k$-good for any $k$.
\end{corollary}

To prove the corollary, just take any non-zero $\alpha$ in $R$, and consider the arithmetic progression $\alpha,2\alpha,3\alpha,4\alpha,\dots$. Then by the previous theorem we obtain that for any $k$, one of these numbers cannot belong to $H_k$. Hence the statement follows.

As always, answering a question opens up several new ones:
\begin{itemize}
\item We have $u(R)=\omega$ or $\infty$ - can we say something about it?
\item Can we say ``how many'' elements of $R$ is the sum of $k$ units?
\item What if we allow linear combinations rather than sums of units?
\item What if units are replaced by elements of bounded norm?
\item What if we require the units to be distinct?
\end{itemize}

First we give some results concerning the first question above.
\begin{itemize}
\item {Belcher (\cite{MR0409409}; 1975), and Ashrafi and V\'amos (\cite{MR2124574}; 2005):} Full characterization of the case $u(R)=\omega$, where $R$ is the ring of integers of ${\mathbb Q}(\sqrt{d})$.
\item {Tichy and Ziegler (\cite{MR2308826}; 2007):} Full characterization of the case $u(R)=\omega$, where $R$ is the ring of integers of ${\mathbb Q}(\sqrt[3]{d})$.
\item {Filipin, Tichy and Ziegler (\cite{MR2490092}; 2008):} Full characterization of the case $u(R)=\omega$, where $R$ is the ring of integers of ${\mathbb Q}(\sqrt[4]{d})$ $(d<0)$.
\item {Frei (\cite{MR2881334}; 2012):} For any number field $K$, there exists a number field $L$ containing $K$, such that the ring of integers of $L$ is generated by its units. This result answers a question of {Jarden and Narkiewicz \cite{MR2309537}}.
\end{itemize}

Now about the question that ``how many" algebraic integers can be obtained as a sum of $k$ units. To make this question precise, we need to introduce some notation. Let $K$ be an algebraic number field with ring of integers $R$, $k$ a positive integer and $x$ a positive real number. Write $u_K(k,x)$ for the number of all classes $[\alpha]$ of associated elements $\alpha\in R$ with $|N(\alpha)|\leq x$ such that $\alpha$ can be written as \[\alpha=\sum\limits_{i=1}^k \varepsilon_i,\] where the $\varepsilon_i$ are units of $R$ and no subsum on the right hand side above vanishes. Now we can give some results concerning $u_K(k,x)$:
\begin{itemize}
\item {Fuchs, Tichy and Ziegler (\cite{MR2540792}; 2009):} Asymptotically precise expression for $u_K(k,x)$.
\item {Frei, Tichy and Ziegler (\cite{MR3260868}; 2014):} A similar result for the case when representations of the form $\alpha=\alpha_1+\cdots+\alpha_k$ are considered, where the $\alpha_i$ are $S$-integers of bounded norm.\footnote{Actually, there was a mistake in a result that was used in \cite{MR2540792} (this is contained in \cite[Lemma 2]{MR2457262}); this gap was successfully closed by Frei, Tichy and Ziegler.}
\end{itemize}

We also mention that there are several results  about the representation of elements of $R$ as linear combinations (rather than sums) of units of $K$. One can see e.g. papers of
Thuswaldner and Ziegler (\cite{MR2825236}; 2011), and Dombek, Hajdu and Peth\H{o} (\cite{MR3217552}; 2014). However, at this point we do not go into details.

Finally, we consider the variant of the problem where in the representations the units have to be \emph{distinct}. If every element of $R$ can be obtained as a sum of distinct units of $K$, then we say that $K$ is a {\sl distinct unit generated} (DUG) field. In this direction we mention the following results:
\begin{itemize}
\item {Jacobson (\cite{MR0160746}; 1964), \'Sliwa (\cite{MR0345929}; 1974)}: ${\mathbb Q}(\sqrt{2})$ and ${\mathbb Q}(\sqrt{5})$ are the only quadratic DUG number fields.
\item {\'Sliwa (\cite{MR0345929}; 1974)}: there are no DUG fields of the form ${\mathbb Q}(\sqrt[3]{d})$; {Belcher (\cite{MR0409409}; 1974)}: there are only seven complex cubic DUG number fields.
\item {Hajdu and Ziegler (\cite{MR3167469}; 2014)}: a finite list of all complex quartic DUG fields. They extended methods of {\'Sliwa, Belcher, Thuswaldner and Ziegler}.
\item {Dombek, Mas\'akov\'a and Ziegler (\cite{MR3283182}; 2015):} extending the above methods further, could make the results of Hajdu and Ziegler much more precise.
\end{itemize}

\section{Conclusion}

We hope that we could show that our groups have worked and work together in a very fruitful and efficient way. We find the results very impressive and it is certainly our wish that the collaboration continues this way. We mention some projects which are currently considered:
\begin{itemize}
\item Diophantine tuples with values in recurrences,
\item polynomial-exponential Diophantine equations,
\item separated-variables Diophantine equations,
\item questions on quadrinomials,
\item lacunary composite polynomials,
\item discrete tomography,
\item and, fortunately, in fact there are many more!
\end{itemize}
The authors would like to thank A. Peth\H{o} and R. Tichy for giving the motivation to such a paper and for their encouragement to turn the presentation given in Gy\H{o}r (Hungary) on August 26, 2015 into a paper. Moreover, we thank the J\'anos Bolyai Mathematical Society, in particular to its President Gyula Katona, for the support in turning this into reality.

The authors are also grateful to the referees for their helpful suggestions.

\end{document}